\documentclass[journal]{IEEEtran}

\usepackage[hyphens]{url}
\usepackage{amsmath, amstext, amssymb}
\usepackage{graphicx}
\usepackage{psfrag}
\usepackage{multirow}
\usepackage{hhline}
\usepackage{eurosym}
\usepackage{color}

\usepackage{array}
\usepackage[latin1]{inputenc}

\begin{document}
%
\title{Security Constrained Optimal Power Flow \\ with Distributionally Robust Chance Constraints}

\author{Line~Roald,~\IEEEmembership{Student Member,~IEEE,}
        Frauke~Oldewurtel,~\IEEEmembership{Member,~IEEE,}
        Bart~Van~Parys,~\IEEEmembership{Student Member,~IEEE,}
        and~Göran~Andersson,~\IEEEmembership{Fellow,~IEEE}

\thanks{L. Roald, F. Oldewurtel and G. Andersson are with the Power Systems Laboratory at the Department
of Electrical Engineering, ETH Zurich, Switzerland. Email: \{roald $\mid$ oldewurtel $\mid$ andersson\}@eeh.ee.ethz.ch.
B. Van Parys is with the Automatic Control Laboratory at the Department
of Electrical Engineering, ETH Zurich, Switzerland. Email: bartvan@control.ee.ethz.ch}}

%



\maketitle

\begin{abstract}
The growing amount of fluctuating renewable in-feeds and market liberalization increases uncertainty in power system operation. To capture the influence of fluctuations in operational planning, we model the forecast errors of the uncertain in-feeds as random variables and formulate a security constrained optimal power flow using chance constraints. The chance constraints limit the probability of violations of technical constraints, such as generation and transmission limits, but require a tractable reformulation. In this paper, we discuss different analytical reformulations of the chance constraints, based on a given set of assumptions concerning the forecast error distributions.
In particular, we discuss reformulations that do not assume a normal distribution, and admit an analytical reformulation given only a mean vector and covariance matrix. 
We illustrate our method with a case study of the IEEE 118 bus system, based on real data from the European system. The different reformulations are compared in terms of both achieved empirical violation probability and operational cost, which allows us to provide a suggestion for the most appropriate reformulation in an optimal power flow setting. For a large number of uncertainty sources, it is observed that the distributions of the line flows and generator outputs can be close to normal, even though the power injections are not normally distributed.
\end{abstract}

\begin{IEEEkeywords}
Renewable integration, Chance Constrained Optimal Power Flow, N-1 security
\end{IEEEkeywords}

%
\IEEEpeerreviewmaketitle

\section{Introduction}
A fundamental tool in power system analysis is the \emph{optimal power flow} (OPF) \cite{stott}.
Several tasks central to power system operation, such as unit commitment, reserve procurement, market clearing and security assessment rely on the solution of an OPF.
The main goal of the OPF is to minimize operational cost, while ensuring secure operation that respects technical limits of the power system.
In current operational schemes, the system is considered secure if it remains within the operational limits during normal operation and during \emph{outage of any single component}. This principle is referred to as the $N-1$ criterion, and is reflected in the OPF through additional constraints, leading to a \emph{security constrained optimal power flow} (SCOPF).
While the $N-1$ criterion secures the system against individual outages, \emph{forecast uncertainty} is another kind of disturbance affecting the system.
Forecast uncertainty arises from unforeseen fluctuations in the power injections, such as inaccurate predictions of load or renewable in-feeds, as well as from short-term electricity trading.
While load profiles are relatively predictable, higher shares of electricity production from renewable sources and liberalization of energy markets (particularly in Europe) have increased the forecast uncertainty by orders of magnitude \cite{zong}.
In current operational planning, uncertainty is usually ignored and uncertain quantities are typically replaced by a forecast value.
While this approach has provided good solutions in the past, the increased levels of uncertainty lead to frequent N-1 violations in real-time operation.
To mitigate these problems, it is proposed to explicitly account for uncertainty during operational planning, in particular while solving the OPF.

There are different approaches to account for uncertainty within the OPF. Robust and worst-case methods, e.g. \cite{joe}, ensure secure operations for all possible forecast errors, but often provide very conservative and thus costly solutions. Stochastic programming methods give the operator more freedom to trade-off cost and security. One example is two stage stochastic programming for unit commitment and reserve scheduling, e.g. \cite{tony}, which minimizes the expected cost of operations based on a set of scenarios. Another example is chance constrained programming, which explicitly limits the probability of constraint violations \cite{maria, chertkov, line}. Since the main goal in short-term operational planning is to ensure secure operations, we consider the latter method, and formulate the OPF as a \emph{probabilistic SCOPF} (pSCOPF) with chance constraints. The acceptable violation probability, which is treated as a design parameter in the optimization problem, allows the operator to choose an appropriate trade-off between cost and security of operations.

Although the pSCOPF allow us to account for uncertainty in a comprehensive way, it is generally hard to reformulate chance constraints as tractable constraints. Two main approaches for reformulation have been applied to the OPF problem, based either on sampling or analytical reformulation.
In \cite{maria}, the SCOPF is formulated as a \emph{joint} chance constrained problem (limiting the probability that any of the constraints are violated), which is reformulated using the scenario approach based on \cite{Campi2006}. The formulation was extended to include market clearing with co-optimization of energy and reserves in \cite{maria2}, where a different sampling based reformulation based on \cite{kostas} was used. Both sampling based reformulations require no knowledge about the underlying distribution, except for availability of a given number of samples (which increases with the problem size).

In contrast, the SCOPF formulated in \cite{line} uses \emph{separate} chance constraints (limiting the probability for each constraint separately) to formulate the pSCOPF. Assuming that the random variables follow a Gaussian distribution, an exact analytical reformulation is obtained. The same type of Gaussian reformulation is performed for an OPF without security constraints in \cite{chertkov}.

While the assumption of a Gaussian distribution limits the applicability of the analytical reformulation from \cite{line}, \cite{chertkov}, the analytical reformulation has some attractive properties.
First, it is scalable to a large number of random variables, as the number of random variables does not influence the problem size or complexity of the OPF itself.
Second, the solution is more transparent than a sample based solution since it is possible to trace the influence of each random variable through the analytical relations.
Finally, the solution based on the analytical reformulation is deterministic, i.e., the OPF will always find the same optimal solution with the same optimal cost. While this might seem trivial, the OPF solution based on the scenario approach is actually random, since it depends on the choice of the samples. The same problem might thus lead to different solutions with different costs, depending on which samples were chosen.

This paper investigates how the good qualities of the analytical reformulation can be preserved, while moving away from the limiting assumption of a Gaussian distribution.
Using optimal probability inequalities, we obtain \emph{distributionally robust} reformulations of the chance constraints. This approach is well-known in operations research and control theory and has been investigated in, e.g.,  \cite{calafiore, popescu, vanparys2013distributionally}.
In \cite{Tyler}, the application of some distributionally robust reformulations to the optimal power flow problem were also discussed. 
Here, we introduce reformulations based on assumptions like unimodality and symmetry of the forecast errors, and explain why those are relevant in the optimal power flow context. We aim to provide 
recommendations for the most suitable reformulations, depending on the sources of uncertainty (e.g., RES fluctuations, load variations, short term trading) and the time frame (e.g., day-ahead planning, real-time operation).
To compare the different analytical reformulations, we introduce the concept of an \emph{uncertainty margin}. The uncertainty margin has a physical interpretation as a security margin against forecast errors, and represents a reduction of available transmission and generation capacity.
A larger uncertainty margin thus increases security, but also the operational cost.
The empirical performance of the proposed reformulations is assessed through a case study for the IEEE 118 bus system, with uncertainties represented through historical forecast errors from the Austrian Power Grid. We investigate which reformulation is the most appropriate for a chance constrained SCOPF problem, considering empirical violation probability, nominal operational cost and the accuracy of the distributional assumptions.

The remainder of this paper is organized as follows: Section II and III present the uncertainty modeling and the formulation of the chance constrained SCOPF. Section IV discusses different analytical reformulations for these chance constraints, and Section V applies them to the power flow equations. Section V demonstrates the proposed formulation in a case study for the IEEE 118 bus system. Section VI summarizes and concludes the paper.

\section{Modeling forecast uncertainty}
Forecast errors arise from uncontrolled power in-feeds that deviate from their forecasted values, caused by, e.g., fluctuations in load or renewable production. The characteristics of forecast errors $\delta_R$ differ between systems, depending on their generation mix, load characteristics and market structure, and also depend on parameters such as the time of the day, or the forecast horizon.
For example, the minute-to-minute variation in wind in-feeds in the central European system follows a Student t-distribution \cite{ENTSOE}, while the day-ahead distribution of wind forecast depends on the forecasted wind power, and is typically non-symmetric. 
The uncertain in-feeds can be modeled as random variables with continuous probability distributions. We define the vector of $n$ uncertain in-feeds as
\begin{equation}
 \tilde{P}_{R} = P_{R}+\delta_R~.
 \label{eq:res}
\end{equation}
Here, the uncertain in-feed $\tilde{P}_R$ is the sum of the forecasted value $P_R$ and a random deviation $\delta_R$. Since the sources of uncertainty differ both within and between systems, the full distribution of $\delta_R$ is generally not known. However, we will assume some partial information about the distribution.
In particular, we assume that the mean $\mu_R \in \mathbb{R}^{n}$ and covariance $\Sigma_R \in \mathbb{R}^{n\times n}$ of the forecast errors exist, and can be estimated either based on historical data or through forecasting methods. We allow for non-zero mean, since forecasts are not necessarily based on the expectation of $P_{R}$, but rather on the most probable realization (which are not the same, e.g., for skewed distributions).

\section{Optimal Power Flow Formulation}
We now introduce the mathematical formulation of the pSCOPF for a system with $n$ buses, based on the formulation in \cite{line}. The sets $\mathcal{G}$, $\mathcal{D}$ and $\mathcal{R}$ represent the conventional generators, the fixed loads and the uncertain in-feeds (consisting of, e.g., in-feeds from wind and solar power plants), respectively. To simplify notation, we assume that there is one generator, load and generator connected at each bus, such that $|\mathcal{G}|=|\mathcal{D}|=|\mathcal{R}|=n$. This assumption is however not necessary for the method itself. The set of transmission lines is denoted by $\mathcal{L}$, and there are $|\mathcal{L}|=n_L$ lines in the system. The contingencies considered for the N-1 security criterion include outage of any line or generator, in total $n_C = n_L + n$ outages.
\subsection{Generator modeling}
The nominal generation output of the generators, $P_G \in \mathbb{R}^{n}$ are the optimization variables of the problem.
In addition to keeping the system balanced in nominal operation conditions, any power deviation arising from either forecast errors or generation outages must be balanced by the generators.
The contribution of balancing energy from each generator can be chosen in different ways. 
Here, we assume that each generator contributes according to its maximum nominal output, similar to \cite{spyros}. When all generators operate, the balancing contribution of each generator $g$ is given by
\begin{equation}
d_{(g)}^i=\frac{P_{G,g}^{max}}{\sum_{j=1}^{n}{P_{G,j}^{max}}}~, 
\label{eq:dg}
\end{equation}
where the superscript $0$ refers to normal operating condition ($i=0$), or situation with line outages ($i \in \mathcal{L}$).
During the outage of generator $i$, the compensation vector of the generators is given by 
\begin{align}
&d_{(g)}^i=\frac{P_{G,g}^{max}}{\sum_{j=1, j\neq i}^{n}{P_{G,j}^{max}}}~ \forall_{\mathcal{G}\setminus i},
&d_{(i)}^i=0~. \label{eq:dg}
\end{align}
The vectors $d^i \in \mathbb{R}^{n}$ thus describe the compensation of any power mismatch in the system, for any outage situation $i$.
We note that by definition of $d$, the system remains balanced after any fluctuations or generator outages.
\subsection{Line flow modeling}
Similar to the setup in \cite{maria}, the line flows are expressed as linear functions of the active power injections in both normal and outage conditions.
\begin{equation}
P_{l}^i=A^i P_{inj}^i, ~~\text{for all }i=0,...,n_C~.
\label{eq:lineflows}
\end{equation}
Here, $A^i\in\mathbb{R}^{n_L\times n}$ describes the relation between the active power injections $P_{inj}^i\in\mathbb{R}^{n}$ and the line flows $P_l^i$ after outage $i$, with $i=0$ being the normal operation condition. $A^i$ is given by
\begin{equation}
A^i = B_f^i \begin{bmatrix} (\widetilde{B}_{bus}^i)^{-1}   ~~~ \bold{0}  \\ ~~\bold{0}  ~~~~~~~~~ 0\end{bmatrix}
\end{equation}
where ${B}_{f}^i\in\mathbb{R}^{n_L\times n}$ is the line susceptance matrix and $\widetilde{B}_{bus}^i\in\mathbb{R}^{n-1\times n-1}$ the bus susceptance matrix (without the last column and row) after outage $i$ \cite{maria}.
The power injections are given by
\begin{align}
P_{inj}^i&=P_G + d^i (P_{G(i)} -\mathbf{1}_{\delta}\delta_R)+P_R + \delta_R - P_D~.
\end{align}
Here, $P_D\in \mathbb{R}^{n}$ is the vector of loads, and the power mismatch due to generation outage $P_{G(i)}$ is non-zero only for $i \in \mathcal{G}$. The vector $\mathbf{1}_{\delta} \in \mathbb{R}^{1\times n}$ is a vector of ones, such that $\mathbf{1}_{\delta}\delta_R$ represents the sum of the forecast errors.

\subsection{Chance constrained optimal power flow}
Using the modeling assumptions presented above, we can formulate the pSCOPF as
\begin{equation}
\min_{P_G} ~c^T P_G
\label{eq:cost}
\end{equation}
subject to
\begin{align}
&\mathbf{1}_{1\times n}(P_G + P_{R} - P_L) = 0  \label{eq:power_balance} \\
&\mathbb{P} [ P_{G(g)} + d^i_{(g)} (P_{G(i)} -\mathbf{1}_{\delta}\delta_R) \leq P_{G(g)}^{max}  ] \geq 1-\varepsilon~, \label{eq:generation upper} \\
&\mathbb{P} [ P_{G(g)} + d^i_{(g)} (P_{G(i)} -\mathbf{1}_{\delta}\delta_R) \geq P_{G(g)}^{min}  ] \geq 1-\varepsilon~, \label{eq:generation lower} \\
&\mathbb{P} [ A^i_{(l,\cdot)}P_{inj}^i \leq  P_{L(l)}^{max}] \geq 1-\varepsilon~, \label{eq:line flow upper} \\
&\mathbb{P} [ A^i_{(l,\cdot)}P_{inj}^i \geq -P_{L(l)}^{max}] \geq 1-\varepsilon~, \label{eq:line flow lower} \\
& \textrm{for } g=1,...,n,~l=1,...,n_L,~i=1,...,n_C~. \nonumber
\end{align}
The objective \eqref{eq:cost} is to minimize generation cost, with $c$ representing the bids of the generators. Constraint \eqref{eq:power_balance} ensures power balance in the system.
The constraints \eqref{eq:generation upper}-\eqref{eq:line flow lower} are the generation and transmission constraints, with $P_G^{min}$ and $P_G^{max}$ being the minimum and maximum generation levels and $P_L^{max}$ being the transmission capacity of the lines.
Those constraints depend on the realization of the random variable $\delta_R$, and are formulated as single chance constraints. The chance constraint ensures that probability of a constraint violation (e.g., a line flow exceeding the limit) remains smaller than $\varepsilon$. We will refer to $\varepsilon$ as the violation probability and to $1-\varepsilon$ as the security level. The value of $\varepsilon$ is an input parameter to the optimization.

\section{Chance constraint reformulation}
\label{sec:reformulations}
To obtain a tractable optimization problem, the chance constraints \eqref{eq:generation upper}-\eqref{eq:line flow lower} must admit a deterministic and tractable reformulation.
These constraints \eqref{eq:generation upper}-\eqref{eq:line flow lower} are all univariate or single chance constraints of the general form
\begin{equation}
\mathbb{P} [ a(P_G) + b(P_G)\delta_R \leq c] \geq 1-\varepsilon~. \label{eq:generic}
\end{equation}
where $a(P_G) \in \mathbb{R}$ and $b(P_G)\in \mathbb{R}^{1\times n}$ are affine functions of the decision variables $P_G$ and $c$ is a constant. The term $a(P_G)$ represents the nominal generation output or the nominal line flows (without forecast errors) and $c$ represents the generation or line flow limit. The vector $b(P_G)$ expresses the influence of the forecast errors $\delta_R$ on the respective constraint.
Regardless of the exact expressions for $b(P_G)$, and for any dimension or distribution of the random vector $\delta_R$, the left hand side of the constraint is a scalar random variable $\delta=a(P_G)+b(P_G)\delta_R$ with mean $\mu(P_G)$ and variance $\sigma(P_G)$ given by
\begin{equation}
\mu(P_G) = a(P_G)+b(P_G)\mu_R~, \quad \sigma(P_G) = \parallel b(P_G)\Sigma_R^{1/2}\parallel_2~. \nonumber
\end{equation}
What is of interest when reformulating the constraint \eqref{eq:generic} is not the distribution of the forecast uncertainty $\delta_R$, but the distribution of $\delta$, which represents the variations in line flows or generation outputs. Depending on the system, $\delta$ might follow different distributions. 
We will now present different distributional assumptions for $\delta$ which are relevant in the context of the SCOPF. The applicability of each assumption depends mainly on the source of uncertainty (e.g., load, renewables or short-term trading), the time frame of the forecast (e.g., day-ahead planning or close to real-time operation) and the availability of data (e.g., historical forecast errors or probabilistic forecasts).
\subsubsection{Normal distribution ($\Phi$)}
The normal distribution is a good distribution model in two different cases. First, when $\delta_R$ follows a multivariate normal distribution (which might be the case, e.g., for load uncertainty), which means $\delta$ will be normally distributed as well. Second, when the number of uncertainty sources is large and not highly correlated, arguments similar to the central limit theorem (e.g., \cite{dasgupta}) imply that the distribution of $\delta$ (which is a weighted sum of $\delta_R$) is expected to be close to a normal distribution.
\subsubsection{Student's $t$-distribution ({$t$})}
When the forecast fluctuations are heavy tailed (e.g., as for the minute-to-minute variability in the European grid \cite{ENTSOE}), the Student's $t$-distribution can be a more appropriate representation. Particularly when considering small violation probabilities ($\varepsilon<0.03$), Student's t distribution provides additional robustness compared to the normal distribution.

In many cases,  only limited knowledge about the distribution of $\delta$ is available. It might therefore be desirable to only assume some general properties of the distribution of $\delta$, rather than a specific distribution. This leads to the following \emph{distributionally robust} reformulations, that are valid for all probability distributions that share the general properties:
\subsubsection{Symmetric, unimodal distributions ($S$)}
If the distribution is likely to be close to normal, but we do not know how close, we can resort to the general assumption of unimodal, symmetric distribution with known mean and covariance.
\subsubsection{Unimodal distributions ($U$)}
In systems where the forecast uncertainty is related mainly to load, wind and PV production, the distribution of $\delta_R$ is likely to be unimodal, with fluctuations centered around the forecasted value.
Under such conditions, it is highly probable that the distribution of $\delta$ is also unimodal. 
\subsubsection{Known mean and covariance ($C$)}
In systems where intra-day electricity trading is not controlled by the transmission system operator, for example in Europe, intra-day transactions introduce uncertainty in the power injections from conventional power plants. The transactions might follow almost any probability distribution, and can even be discrete. In this case, we reformulate the chance constraint based only on a known (and finite) mean and covariance.

For all distributional assumptions 1) - 5), the chance constraint \eqref{eq:generic} can be reformulated to the following analytic expression 
\begin{equation}
a(P_G) \leq c - b(P_G)\mu_R - f^{-1}(1-\varepsilon)\parallel b(P_G)\Sigma_R^{1/2}\parallel_2~. \label{eq:reformulated}
\end{equation}
Analyzing \eqref{eq:reformulated}, we see that the left part, $a(P_G)\leq c$, represents the ``nominal'' constraint, i.e., the constraint we would obtain if we neglect the forecast uncertainty. The second and third term represents a reduction of the nominally available capacity $c$, which is necessary to secure the system against forecast deviations. This reduction can thus be interpreted as a security margin against uncertainty, i.e., an \emph{uncertainty margin}. Notice that the larger $f^{-1}(1-\varepsilon)$, the larger the uncertainty margin.

Depending on which assumption 1) - 5) is deemed appropriate, we define $f^{-1}(1-\varepsilon)$ according to either an inverse cumulative distribution function (for known distributions 1), 2)) or a probability inequality (when only partial information is available 3) - 5)). The exact expressions for $f^{-1}(1-\varepsilon)$ are shown in Table \ref{tableI}, and their derivations as well as the derivation of \eqref{eq:reformulated} are given in the Appendix.
We note that for 1) and 2), the reformulation is tight (the chance constraint holds with equality).
The distributionally robust reformulations 3)-5) are typically not tight, and will usually lead to empirical violation probabilities lower than $\varepsilon$.

{\setlength{\extrarowheight}{6pt}
\setlength\belowcaptionskip{-30pt}
\begin{table} [h!]
\caption{Expressions for $f^{-1}(1-\varepsilon)$.}
$\Phi$: Cumulative distribution function of the standard normal distribution. $t_{\nu,\sigma_T}$: Cumulative distribution function of the Student t distribution with zero mean, $\nu$ degrees of freedom and scale parameter $\textstyle{\sigma_T = (\nu-2)/\nu}$.
\label{tableI}
\centering
\begin{tabular}{|l|l|}
\hline
{1) Normal }       & $f_{\Phi}^{-1}(1-\varepsilon) = \Phi^{-1}(1-\varepsilon)$ \\[2pt]
\hline
{2) Student's $t$}            & $f_{t}^{-1}(1-\varepsilon)    = t_{\nu, \sigma_T}^{-1}(1-\varepsilon)$ \\[2pt]
\hline
{3) Symmetric, unimodal}        & $f_{S}(1-\varepsilon)    = \begin{cases} \sqrt{\frac{2}{9\varepsilon}}  ~ &\text{for } 0\leq\varepsilon\leq\frac{1}{6} \\
                                                                            \sqrt{3}(1-2\varepsilon)  ~ &\text{for } \frac{1}{6}<\varepsilon<\frac{1}{2} \\
                                                                            0 ~ &\text{for } \frac{1}{2}\leq\varepsilon\leq 1 \end{cases}$ \\[6pt]

\hline
{4) Unimodal}                  & $f_{U}^{-1}(1-\varepsilon)    = \begin{cases} \sqrt{\frac{4}{9\varepsilon}-1}  \quad &\text{for } 0\leq\varepsilon\leq\frac{1}{6} \\
                                                                            \sqrt{\frac{3(1-\varepsilon)}{1+3\varepsilon}} \quad &\text{for } \frac{1}{6}<\varepsilon\leq 1 \end{cases}$ \\[6pt]
\hline
{5) Mean, covariance}                 & $f_{C}^{-1}(1-\varepsilon)    = \sqrt{\frac{1-\varepsilon}{\varepsilon}} \quad \text{for } 0\leq \varepsilon\leq 1$\\[3pt]
\hline
\end{tabular}
\end{table}}

\begin{figure}
\includegraphics[width=0.95\columnwidth]{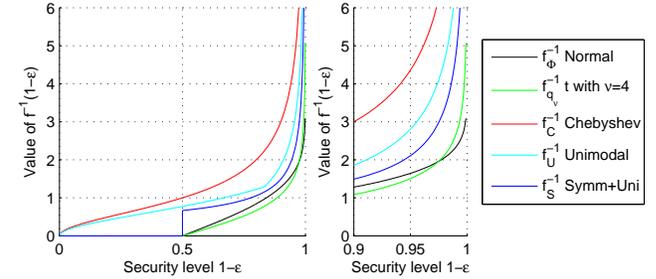}
\centering
\caption{Values of $f^{-1}(1-\varepsilon)$ for the normal distribution, the Student t distribution with 4 degrees of freedom, the Chebyshev inequality, unimodal distributions and symmetric, unimodal distributions. The left part shows all security levels, while the right part is a zoom in on high security levels. }
\label{fig_finv}
\end{figure}

Since the reformulations 1) - 5) differ only in the definition of $f^{-1}(1-\varepsilon)$, we can compare them by comparing the value of $f^{-1}(1-\varepsilon)$ for different $\varepsilon$.
In Fig. \ref{fig_finv}, $f^{-1}(1-\varepsilon)$ is plotted against the security level $1-\varepsilon$.
We observe that all $f^{-1}(1-\varepsilon)$ increase as $\varepsilon$ decreases, indicating that a larger uncertainty margin is necessary to achieve a lower violation probability.
With more information, tighter probabilistic bounds can be defined and thus a lower value of $f^{-1}(1-\varepsilon)$ is necessary to ensure the desired security level (i.e., $f_C^{-1}>f_U^{-1}>f_S^{-1}$). The lowest values are obtained when we assume knowledge of the actual distribution, i.e., for the normal and the Student's $t$ distribution. Note that all reformulations assuming symmetry have $f_S(0.5)=f_\Phi(0.5)=f_t(0.5)=0$.
Finally, Student's $t$ distribution has a more pronounced peak and heavier tails than the normal distribution. This is reflected in that for lower security levels, $f_{t}^{-1}(1-\varepsilon)<f_{\Phi}^{-1}(1-\varepsilon)$, while at high security levels, $f_{t}^{-1}(1-\varepsilon)>f_{\Phi}^{-1}(1-\varepsilon)$.

\section{Reformulated constraints}
With the reformulation presented above, we can reformulate the chance constraints \eqref{eq:generation upper}-\eqref{eq:line flow lower} as
\begin{align}
&P_{G(g)} + d^i_{(g)} P_{G(i)}  \leq P_{G(g)}^{max}- \label{eq:generation upper_ref} \\
& \quad d^i_{(g)}\mathbf{1}_{\delta}\mu_R - f^{-1}(1-\varepsilon) \parallel d^i_{(g)}\mathbf{1}_{\delta}\Sigma_R^{1/2}\parallel_2, \nonumber \\
&P_{G(g)} + d^i_{(g)} P_{G(i)}  \geq P_{G(g)}^{min}- \label{eq:generation lower_ref} \\
& \quad d^i_{(g)}\mathbf{1}_{\delta}\mu_R + f^{-1}(1-\varepsilon) \parallel d^i_{(g)}\mathbf{1}_{\delta}\Sigma_R^{1/2}\parallel_2, \nonumber \\
&A^i_{(l,\cdot)}(P_G + d^i P_{G(i)} + P_R - P_D) \leq  P_{L(l)}^{max}- \label{eq:line flow upper_ref}\\
& \quad A^i_{(l,\cdot)}(I\!-\!d^i \mathbf{1}_{\delta})\mu_R  - f^{-1}(1\!-\!\varepsilon)\!\parallel \! A^i_{(l,\cdot)}(I\!-\!d^i \mathbf{1}_{\delta})\Sigma_R^{1/2}\!\parallel_2, \nonumber \\
&A^i_{(l,\cdot)}(P_G + d^i P_{G(i)}+P_R - P_D) \geq  -P_{L(l)}^{max}- \label{eq:line flow upper_ref}\\
& \quad A^i_{(l,\cdot)}(I\!-\!d^i \mathbf{1}_{\delta})\mu_R  + f^{-1}(1\!-\!\varepsilon)\!\parallel \! A^i_{(l,\cdot)}(I\!-\!d^i \mathbf{1}_{\delta})\Sigma_R^{1/2}\!\parallel_2, \nonumber \\
& \textrm{for } g=1,...,n,~l=1,...,n,~i=1,...,n_C~, \nonumber
\end{align}
where $I\in\mathbb{R}^{n\times n}$ is the identity matrix.
When comparing \eqref{eq:generation upper_ref} - \eqref{eq:line flow upper_ref} to \eqref{eq:reformulated}, we recognize the same structure.
The first part represent the constraint we would obtain by neglecting the forecast uncertainty, where as the second and third term on the right hand side represent the uncertainty margin.

Since a higher uncertainty margin leads to a reduction in the available transmission and generation capacity, a higher uncertainty margin will not only reduce the probability of violation, but also increase the nominal cost of operation (i.e., the cost of the pSCOPF). The acceptable violation probability $\varepsilon$ and the distributional assumption (which defines the function $f^{-1}(1-\varepsilon)$) should therefore be chosen carefully to obtain a good trade-off between security against forecast errors and cost of operation.

%

Note that the reformulated chance constraints \eqref{eq:generation upper_ref} - \eqref{eq:line flow upper_ref} are linear, since the uncertainty margin is not dependent on any decision variables and can be pre-computed. The pSCOPF problem \eqref{eq:cost}, \eqref{eq:power_balance}, \eqref{eq:generation upper_ref} - \eqref{eq:line flow upper_ref} is thus a linear program with the same computational complexity as a traditional DC SCOPF.

\section{Case Studies}
\label{sec:case_study}
\setcounter{subsubsection}{0}
The purpose of this case study is to demonstrate the chance-constrained SCOPF, and investigate which distributional assumptions are most appropriate for power systems operation.
We base our study on the IEEE 118-bus system \cite{118busdata}, with a few modifications as follows. The generation cost is assumed to be linear, and is based on the linear cost coefficients of the data provided with Matpower 4.1 \cite{matpower}. Although the formulation could be extended to include unit commitment, it is not considered here.  Therefore, the minimum generation output of the conventional generators is set to zero.
The forecast uncertainty $\delta_R$ is modeled based on historical data for 1 year from the Austrian Power Grid (APG). We define the forecast error as the difference between the the so-called DACF (Day-Ahead Congestion Forecast) and the snapshot (the real-time power injections) for all hours and buses with available data (8492 data points for 28 buses). Since the system is constantly evolving and might exhibit seasonal patterns, we assume that the power system operator only uses data from the past three months. We use two three-month periods to define the forecast uncertainty for this case study, such that we obtain 2207 data samples for a total of 54 buses.

\begin{figure}[!t]
\includegraphics[width=1\columnwidth]{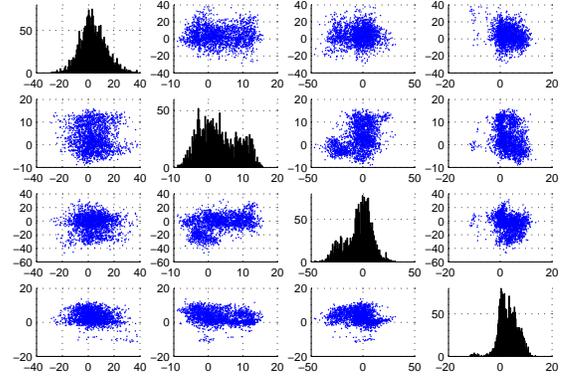}
\centering
\caption{Forecast errors for 4 selected nodes of case study. The diagonal plots show the histograms of the forecast errors (x-axis: deviation in MW, y-axis: number of occurences), while the off-diagonal plots show the scatter plots between two corresponding forecast errors (x- and y-axis: deviation in MW).}
\label{fig_inputs}
\end{figure}

The historical data was assigned to different load buses throughout the system, and modified such that the standard deviation corresponds to 20 \% of the forecasted load. The mean $\mu_R$ and covariance $\Sigma_R$ used in the pSCOPF were calculated based on this modified data (i.e, assuming perfect knowledge of $\mu_R,~\Sigma_R$). Fig. \ref{fig_inputs} shows the forecast errors from some representative nodes, including the histograms and pair-wise scatter plots of the forecast errors. By inspection, it is clear that the forecast errors are not normally distributed.

In the following, we assess how the different distributional assumptions impact the solution of the pSCOPF. We solve the pSCOPF for all five reformulations assuming an acceptable violation probability of $\epsilon = 0.1$. The results are compared with each other and to the solution of the corresponding deterministic SCOPF. To assess the quality of the solution, we compare the number of empirical constraint violations (based on the historical samples) and the relative cost of the solutions. Further, we run statistical tests to check if the data is normally or unimodally distributed, and investigate the accuracy of the estimated uncertainty margins.

\subsubsection{Number of empirical violations}
\begin{figure}[!t]
\includegraphics[width=0.9\columnwidth]{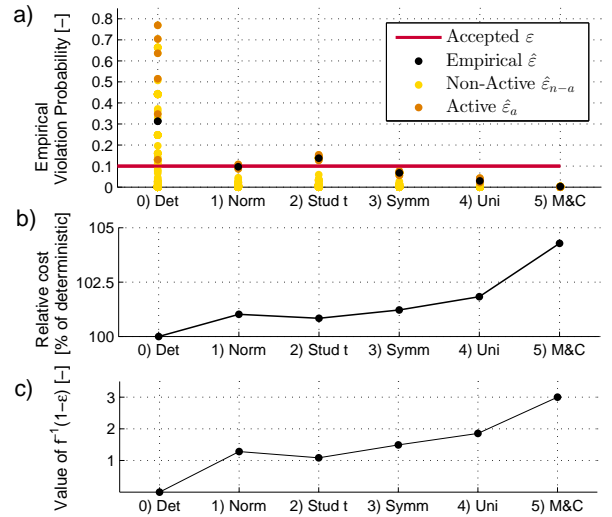}
\centering
\caption{Results derived from the different SCOPF solutions. From left to right: 0) the deterministic SCOPF, and the pSCOPF based on the assumption of 1) a normal distribution, 2) a Student t distribution, 3) a symmetric, unimodal distribution, 4) a unimodal distribution and 5) known mean and covariance. From top to bottom, the figure shows a) empirical violation probabilities for the all constraints, b) nominal dispatch cost, and c) value of the function $f^{-1}(1-\varepsilon)$ for $\varepsilon = 0.1$.}
\label{cost}
\end{figure}

The empirical violation probabilities $\hat\varepsilon$ are evaluated for all constraints based on the 2207 data samples. The results are shown in Fig. \ref{cost} a), starting with 0) the deterministic solution, then the solution based on 1) a normal and 2) a Student t distribution, and then the distributionally robust solutions 3)-5). From left to right, we thus assume lessened knowledge about the distribution. The violation probabilities of the non-active constraints $\hat\varepsilon_{n-a}$ and active $\hat\varepsilon_{a}$ constraints are plotted in yellow and orange, respectively. The average violation probability $\hat\varepsilon_{avg}$ of the active constraints is plotted in black.

As seen on the left, the empirical violation probability of the deterministic solution is very high, with some constraints violating 80\% of the cases. This highlights the need for probabilistic methods to avoid frequent violations of operational limits.
The pSCOPF solutions have much lower violation probabilities. The active constraints (where the distance from the nominal flow to the flow limit is given by the uncertainty margin) have a higher empirical violation probability than the non-active constraints (which have some additional margin). The solution based on a normal distribution violates the accepted violation probability $\varepsilon < 0.1$ for some constraints, but the violation is small, $\hat\varepsilon<0.11$, and the average violation probability is acceptable $\hat\varepsilon_{avg}<0.1$. The solution based on a Student t distribution, which assumes a more peaked distribution, has larger empirical violation probabilities $\hat\varepsilon > 0.15$. The distributionally robust solutions oversatisfy the accepted violation probability $\varepsilon$, with $\hat\varepsilon_{avg}=0.07$ for the symmetric, unimodal solution, $0.03$ for the unimodal solution and $0.0025$ for the solution based only on mean and covariance. Although the chance constraints are satisfied, it does not necessarily imply that the underlying assumption (e.g., symmetry and unimodality) is accurate. Since the reformulations are distributionally robust, we might get a low empirical violation probability, even if we assumed the wrong family of distributions.

\subsubsection{Operational cost}
Fig \ref{cost} b) shows the generation cost obtained with the pSCOPF solutions, normalized by the cost of the deterministic problem (shown to the left).
All probabilistic solutions have higher cost than the deterministic solution, showing that the consideration of uncertainties increase the nominal cost of operation. The reformulations which assume more knowledge about the distribution 1), 2) lead to lower cost than the more general reformulations 3)-5). The most expensive solution is obtained for reformulation 5), which only assumes knowledge of mean and covariance.

The cost differences are explained by the different values $f^{-1}(1-\varepsilon)$, which defines the uncertainty margin and thus the constraint tightening. The value of $f^{-1}(1-\varepsilon)$ is plotted in Fig. \ref{cost} c).
Comparing Fig. \ref{cost} a), b) and c), we observe how a larger $f^{-1}(1-\varepsilon)$ leads to an increase in nominal operation cost, but at the same time reduces the empirical violation probability. This highlights two important aspects of the pSCOPF. First, we need to define $\varepsilon$ such that it reflects a reasonable trade-off between cost and security. Second, we want to achieve an empirical violation probability $\hat\varepsilon$ as close as possible to the accepted violation level $\varepsilon$. A reformulation with too many violations ($\hat\varepsilon>>\varepsilon$) leads to unsecure operations, but at the same time a too conservative solution ($\hat\varepsilon<<\varepsilon$) will lead to unnecessary high cost, and possibly infeasibility if we want to ensure a low violation probability (e.g., $\varepsilon\approx 0.01$).


\subsubsection{Testing the distributional assumptions and the accuracy of estimated uncertainty margin}
Since the transmission and generation constraints are enforced as separate chance constraints with deviations $\delta$ (defined as a weighted sum of the random variables $\delta_R$), each constraint has a univariate distribution related to it.
To assess whether or not our assumptions about those distributions are correct, we run statistical tests. In particular, we use the Shapiro-Wilk test \cite{shapiroWilk} to test if the distribution is normal, and Hartigans dip test \cite{hartigan} to test unimodality, using the implementation in R \cite{citeR}.
The test output is a p-value between 0 and 1, which indicates how probable it is that the data comes from a normal or a unimodal distribution, respectively. Typically, the hypothesis (normality or unimodality) is accepted for p-values above $p>0.95$, and rejected for p-values $p<0.05$.
In between, we can neither reject nor confirm the hypothesis.

\begin{figure}[!t]
\includegraphics[width=1\columnwidth]{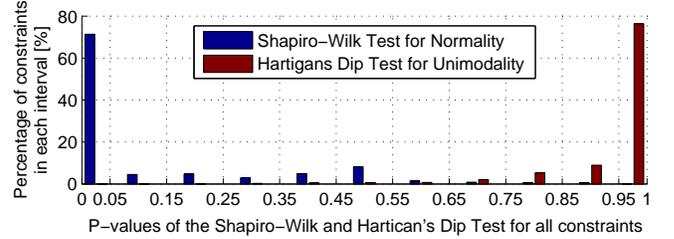}
\centering
\caption{P-values obtained from Hartigans dip test for unimodality and the Shapiro-Wilk test of normality. The histogram show the percentage of all constraints in each 0.05 interval. A high p-value indicate a high probability that the distribution is unimodal or normal, respectively. Based on the results, it seems highly probable that the distributions are unimodal (p-values close to 1), while it is unlikely that the data comes from a normal distribution (p-values close to 0).}
\label{DipTest}
\end{figure}

In Fig. \ref{DipTest}, the p-values from both tests are plotted as a histogram. The bars show the percentage of constraints with p-values in the indicated p-value interval. We observe that the p-values from the Shapiro-Wilk test are below the threshold $p<0.05$ for most constraints, while the p-values from Hartigans Dip Test are above $p>0.95$ for the majority of the constraints. We thus conclude that unimodality is a reasonable assumption, while the original data is probably not normally distributed.

\begin{figure}[!t]
\includegraphics[width=0.9\columnwidth]{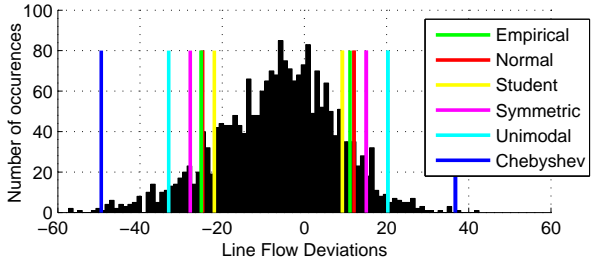}
\centering
\caption{Histogram of post-contingency line flow deviations. The uncertainty margins are computed empirically (green), and for 1) a normal distribution (red), 2) a Student t distribution (yellow), 3) a symmatric, unimodal distribution (magenta), 4) a unimodal distribution (light blue) and 5) a distribution where only the mean and covariance are known (dark blue).}
\label{histMargins}
\end{figure}

Although the statistical test rejects normality, the normal distribution might still be a good assumptions for the parts of the distribution which we are interested in. In Fig. \ref{histMargins}, the empirical distribution of the line flow deviations for one active transmission constraint is shown. This constraint had the lowest p-value among the active constraints in the Shapiro-Wilk test. Fig. \ref{histMargins} also show the empirical uncertainty margin (corresponding to an empirical violation probability $\hat\varepsilon=0.1$), as well as the uncertainty margins obtained with the five distributional assumptions. We observe that the uncertainty margins based on the normal distribution (plotted in red) match very closely to the empirical margins (plotted in green). The Student t distribution underestimates the margin and the distributionally robust reformulations lead to too high margins.

\subsubsection{Discussion}
Based on the above results, we conclude that the reformulation based on a normal distribution can provide a good trade-off between cost and security, particularly in systems with a large number of uncertainty sources.
Although assuming a normal distribution does not guarantee empirical violation probabilities $\hat\varepsilon<\varepsilon$, the assumption might be useful when $\varepsilon$ can be interpreted as a guideline rather than a hard limit. Statistical tests and assessments as in Fig. \ref{histMargins} can be used to assess whether the normal distribution is a reasonable approximation.
If the system operator wants a higher confidence in enforcing the actual violation level $\varepsilon$ and is willing to tolerate a larger increase in operational cost, assuming a unimodal distribution would be reasonable.

\section{Summary and Conclusion}
This paper discusses different analytic reformulations for chance constraints and their applicability in the pSCOPF context. The chance constraints are reformulated either by assuming a known probability distribution (such as normal or Student t distribution) or by using distributionally robust reformulations assuming general properties of the distribution (i.e., known mean and variance, symmetry, unimodality).

The reformulated chance constraints all have a similar form, and are easily comparable. They are similar to the nominal constraints of the deterministic problem, except for the uncertainty margin (a security margin against forecast deviations), which represents a reduction of the transmission or generation capacity. With a larger uncertainty margin, the probability of violations decreases, but the nominal operational cost increases. Therefore, it is desirable to find a reformulation which leads to an uncertainty margin which is sufficiently large, yet as small as possible.


In the case study based on the IEEE 118 bus system and forecast errors from Austria, the trade-off between security and cost is highlighted. Although the choice of reformulation differs between systems with different uncertainty characteristics, we show that the normal distribution might be a good approximation in cases where the acceptable violation probability can be interpreted as a guideline, rather than a hard constraint. If the transmission system operator wants to enforce the violation probability as a strict limit, choosing a more conservative, distributionally robust reformulations based on, e.g., unimodality will provide more confidence.

In general, we believe that the pSCOPF with analytically reformulated chance constraints provides a transparent and scalable approach to assess the effect of uncertainty in power system operational planning. Future work will investigate how the approach can be extended towards corrective control actions for uncertainties (e.g., HVDC and PSTs), to further reduce the cost of handling uncertainty. Further, we plan to investigate extensions towards AC power flow.

\section*{Acknowledgment}
This research work described in this paper has been partially carried out within the scope of the project "Innovative tools for future coordinated and stable operation of the pan-European electricity transmission system (UMBRELLA)", supported under the 7th Framework Programme of the European Union, grant agreement 282775. We thank our partners, especially Austrian Power Grid, for providing historical data.

\appendix

\subsection{Exact tractable reformulations of single chance constraints}

We discuss here how the chance constraint \eqref{eq:generic} can be reformulated as the deterministic constraint \eqref{eq:reformulated}. Using the equality $\delta = a(P_G) + b(P_G) \delta_R$ it is clear that constraint \eqref{eq:generic} is equivalent to
\begin{equation}
\label{eq:alternative1}
	\mathbb P[\delta < c] \geq 1-\epsilon
\end{equation}
where the constraint \eqref{eq:alternative1} should be satisfied for all distributions $\mathbb P$ consistent with the distributional assumptions made on $\delta$. The constraint \eqref{eq:alternative1} can equivalently be represented as
\begin{equation}
\label{eq:alternative2}
	\mathbb P\left[\frac{\delta-\mu(P_G)}{\sigma(P_G)} < \frac{c-\mu(P_G)}{\sigma(P_G)}\right] \geq 1-\epsilon
\end{equation}
where it can be remarked that the scaled random variable $\delta_n:=\left(\delta-\mu(P_G)\right)/\sigma(P_G)$ has zero mean and unit variance by construction.

In order to unify the analysis for all distributional assumptions made in Section \ref{sec:reformulations}, we consider the general situation in which the distribution $\mathbb P$ of $\delta_n$ is merely known to belong to a set of distributions $\mathcal P$. Specifically, we have that $\mathcal P$ correspond to $\{N(0, 1)\}$ in the Gaussian case 1) and to $\{ t_\nu\left(0, \sqrt(\frac{\nu - 2}{\nu})\right) \}$ in the Student's $t$ case 2) with $\nu$ degrees of freedom. In the former two cases, the set $\mathcal P$ is a singleton as the distribution of $\delta_n$ is known. In the most general case 5), the set $\mathcal P$ consists of all distributions of zero mean and unit variance. The set $\mathcal P$ is additionally required to contain only unimodal or unimodal symmetric distributions in the unimodal case 4) and the unimodal symmetric case 3), as both notions are scale invariant. That is, if $\delta$ is unimodal or symmetric unimodal than $\alpha \delta + \beta$ is unimodal or symmetric unimodal as well for any $\alpha$ and $\beta$ real numbers.

Given a set of distributions $\mathcal P$ representing the distributional assumptions made, we will first show that the chance constraint \eqref{eq:generic} is equivalent to the deterministic constraint \eqref{eq:reformulated} for
\(
	f_\mathcal P(k) := \inf_{\mathbb P \in \mathcal P} ~\mathbb P \left [ \delta_n < k \right]
\).

We can trivially rewrite constraint \eqref{eq:alternative2} equivalently as the constraint $\mathbb P[\delta_n < \left(c-\mu(P_G)\right)/\sigma(P_G)] \geq 1- \epsilon$ for all $\mathbb P\in \mathcal P$. Using the definition of $f_{\mathcal P}$ we can demand alternatively that $f_{\mathcal P}(\left(c-\mu(P_G)\right)/\sigma(P_G)) \geq 1- \epsilon$. Finally as the function $f_{\mathcal P}$ is increasing it has a well defined generalized inverse $f_{\mathcal P}^{-1}(\lambda) = \inf \,\{k ~\vert~ f_{\mathcal P}(k)\geq \lambda\}$. From the definition of $f$ it follows that we must have that $\left(c-\mu(P_G)\right)/\sigma(P_G) \geq f_{\mathcal{P}}^{-1}(1-\epsilon)$. After reordering the terms we obtain
\[
	a(P_G) + b(P_G) \mu_R \leq c -f_{\mathcal P}^{-1}(1-\epsilon)\Vert b(P_G)\Sigma_R^{1/2}\Vert_2.
\]

Finally in the remainder of this section, we show for the distributional assumptions discussed above we obtain the results mentioned in Table \ref{tableI}.

\begin{enumerate}
\item[1--2)] For non-atomic distributions $\mathbb P$ it follows by continuity that $f_{\{\mathbb P\}} = \mathbb P[\delta_n < k] = \mathbb P[\delta_n\leq k]$  explaining the first two entries in Table \ref{tableI}.

\item[3)] Let the set $\mathcal P$ corresponds to the set of all symmetric unimodal distributions $S$ with zero mean and unit variance then we can leverage the classical Gauss bound \cite{gauss1821theoria, vanparys2015generalized}.
Indeed, in this situation we have that the mode of the distributions is known as it coincides with the mean because of symmetry. For any $k>0$, we have again by symmetry that the equality $\mathbb P[\delta_n \geq k] = \frac12 \mathbb P[\vert \delta_n \vert \geq k]$ holds. As the mode of $\mathbb P$ is known, the classical Gauss bound \cite{gauss1821theoria, vanparys2015generalized} can be used to establish
\begin{eqnarray}
	f_{U}(k) &=& 1 - \sup_{\mathbb P \in {S}} ~\mathbb P \left [ \delta_n \geq k \right]
    \nonumber \\
    &=& \left\{
	\begin{array}{c l}
	1- \frac12 \sup_{\mathbb P \in {U}} ~\mathbb P \left [ \vert\delta_n\vert \geq k \right] & \mathrm{if~} k> 0 \\
	0 & \mathrm{otherwise}
	\end{array} \right.
	\nonumber \\
    &=&
	\left\{
	\begin{array}{c l}
	\frac{9k^2-2}{9k^2} & \mathrm{if~} k\geq \sqrt{\frac43} \\
	\frac12 +\frac{k}{2\sqrt 3} &  \mathrm{}{k> 0} \\
	0 & \mathrm{otherwise}
	\end{array} \right.
\end{eqnarray}

\item[4)] The unimodal case, in which $\mathcal P = U$ consists of zero mean and unit variance unimodal measures, can be dealt with using the one-sided Vysochanskij--Petunin inequality \cite{vysochanskij1985improvement}, i.e.\
\begin{eqnarray}
	f_{U}(k) &=& 1 - \sup_{\mathbb P \in {U}} ~\mathbb P \left [ \delta_n \geq k \right]
    \nonumber \\
    &=& 1 - \left\{
	\begin{array}{c l}
	\frac{4}{9\left(1+k^2\right)} & \mathrm{if~} k\geq \sqrt{\frac53} \\
	1 - \frac4 3 \frac{k^2}{1+k^2} & \mathrm{if ~} {k\geq 0} \\
	1 & \mathrm{otherwise}
	\end{array} \right.
	\nonumber \\
    &=&
	\left\{
	\begin{array}{c l}
	1 - \frac{4}{9\left(1+k^2\right)} & \mathrm{if~} k\geq \sqrt{\frac53} \\
	\frac4 3 \frac{k^2}{1+k^2} &  \mathrm{}{k\geq 0} \\
	0 & \mathrm{otherwise}
	\end{array} \right.
\end{eqnarray}

\item[5)] Lastly, when the set $\mathcal P$ corresponds to the set $C$ containing all distributions of zero mean and variance, the classical Cantelli inequality \cite{cantelli1910intorno} establishes that
\begin{eqnarray}
	f_{C}(k) &=& 1 - \sup_{\mathbb P \in {C}} ~\mathbb P \left [ \delta_n \geq k \right]
    \nonumber\\
    &=& 1 - \left\{
	\begin{array}{c l}
	\frac{1}{1+k^2} & \mathrm{if~} k\geq 0 \\
	1 & \mathrm{otherwise}
	\end{array} \right.
    \nonumber \\
	&=&
	\left\{
	\begin{array}{c l}
	\frac{k^2}{1+k^2} & \mathrm{if~} k\geq 0 \\
	0 & \mathrm{otherwise}
	\end{array} \right.
\end{eqnarray}
which after taking the inverse gives us the corresponding result in Table \ref{tableI}.
\end{enumerate}
\bibliographystyle{IEEEtran}
\bibliography{20150326_bib_pSCOPF}

%
%
%



\enlargethispage{-5in}


\end{document}